\documentclass[11pt]{article}
\usepackage{graphicx}
\usepackage{amssymb,amsfonts,amsthm}
\usepackage{amsmath,amsthm,amssymb}
\usepackage{amsfonts}

\newtheorem{theorem}{Theorem}[section]
\newtheorem{lemma}[theorem]{Lemma}

\theoremstyle{definition}

\newtheorem{rk}[theorem]{Remark}

\newcounter{ppp}

\newcommand{\la}{\langle}
\newcommand{\ra}{\rangle}

\newcommand{\me}{\medskip}

\newcommand{\tool}{\stackrel{\ell}{\too} }

\newcommand{\eee}{{\cal E}}
\newcommand{\aaa}{{\cal A}}

\newcommand{\iv}{^{-1}}

\newcommand{\too}{\to }

\newcommand{\xxx}{{\cal X} }

\newcommand{\base}{\mathrm{base}}
\begin{document}

\renewcommand{\theequation}{\thesection.\arabic{equation}}

\title{Groups with quadratic-non-quadratic Dehn functions }
 \author{A.Yu. Ol'shanskii \thanks{The author was supported in part by
the NSF grant DMS 0245600 and by the Russian Fund for Basic
Research 02-01-00170.}}

\date{}
\maketitle

\begin{abstract}
We construct a finitely presented group $G$ with non-quadratic
Dehn function $f$ majorizable by a quadratic function on arbitrary
long intervals.

\end{abstract}

\section{Introduction}

Recall that the {\em Dehn function} of a finite presentation $\la
X\mid R\ra$ of a group $G$ is the smallest function $f:{\mathbb
N}\to{\mathbb N}$ such that any word of length at most $n$ in $X$
that represents the identity of $G$ is freely equal to a product of
at most $f(n)$ conjugates of elements of $R$. The Dehn functions
$f_1, f_2$ of any two finite presentations of the same group $G$ are
{\em equivalent}, that is $f_2(n)\le Cf_1(Cn) +Cn+C$,
$f_1(n)<Cf_2(Cn)+Cn+C$ for some constant $C$. As usual, we do not
distinguish equivalent functions. The Dehn function can also be
defined as the smallest {\em isoperimetric} function of the
presentation: that is the smallest function $f(n)$ such that the
area (i.e. the number of 2-cells) of a minimal van Kampen diagrams
over $\la X\mid R\ra$ having perimeter (i.e. the combinatorial
length of the contour) at most $n$ does not exceed $f(n)$. The
connections of the properties of Dehn functions, on the one hand, to
the asymptotic geometry of groups and spaces, and, on the other
hand, to the computational complexity of the algorithmic word
problem,  are discussed in \cite{Gr3}, \cite{BORS} and \cite{OS1}.

The class of increasing functions which, up to equivalence, can be
represented as Dehn functions of groups is vast (see \cite{SBR},
\cite{Bri}, \cite{BrBr}, \cite{OS}), but there is one gap in the
scale of their rates: if a Dehn function of a group $G$ is
subquadratic, then it is linear, and so $G$ is a word hyperbolic
group \cite{Gr2}, \cite{Ol}, \cite{Bo}.

The goal of this paper is to give an example of a group whose Dehn
function is not majorized on ${\mathbb N}$ by a quadratic function
but is smaller than a quadratic function on arbitrary long
intervals of natural numbers.

\begin{theorem} \label{th1}
Let $\Psi(n) = n^2 \log' n /\log'\log'n$ where $n\ge 0$ and
$\log'n=\max (\log_2 n, 1)$. There is a finitely generated group
$G$ whose Dehn function $f(n)$ satisfies the following properties:

(1) $c_1n^2\le f(n)\le c_2 \Psi(n)$ for some positive constants
$c_1, c_2$ and all sufficiently large $n$;

(2) there is a sequence $n_i\to \infty$ with $f(n_i)\ge
c_3\Psi(n_i)$ for a positive $c_3$ and every $n_i$;

(3) there is a sequence $n'_i\to \infty$ with $ f(n'_i)\le c_4
(n'_i)^2$ for a positive $c_4$ and every $n'_i$;

Moreover:

(4) there are sequences of positive numbers $d_i\to\infty$ and
$\lambda_i\to\infty$ such that $f(x)\le c_4 x^2 $ for arbitrary
integer $x\in \cup_{i=1}^{\infty} [\frac{d_i}{\lambda_i} ,\;
\lambda _i d_i]$,

(5) there is a positive constant $c_5$ such that for every $n_i$
defined in (2), and for every integer $x$ with $x\le c_5 n_i$, we
have $f(x)\le c_4 n_i^2$, in particular, $f([c_5 n_i])/f(n_i)\to
0$.

The group $G$ is a multiple HNN extension of a free group.
\end{theorem}

It will be clear from the proof that the sequences $(n_i)$,
$(n'_i)$ and $(d_i)$ have double exponential growth. Notice that
such sequences cannot grow as an ordinary exponential function (or
slower), because any function $f$, up to equivalence, is
determined by its values $f(a^i)$, $i=1,2,\dots$, if $a>1$.

The unusual almost quadratic behavior of the Dehn function, and
especially the properties (4) and (5), lead in \cite{OS2}, to the
solution of a well known problem about asymptotic cones of groups.
Namely, property (4) guarantees that a cone defined by the scaling
sequence $(d_i)$ is simply connected, and property (5) implies
that a cone defined by the sequence $n_i$ is not simply connected.
(See \cite{OS2} for the details.)

To prove Theorem \ref{th1} we construct $G$ as a special multiple
HNN extension of a free groups, namely, an $S$-machine. Starting
with \cite{SBR}, Sapir's S-machines are applied to a number of group
theoretical tasks. In Section 3, we recall the definition and basic
property of an auxiliary adding machine $Z(A)$ introduced in
\cite{OS}. In a sense, the main machine $M$ defined in Section 4, is
composed from various copies of  adding machines.

It is seen from the definition of $M$ that given number $n$, this
machine can produce a computation $W_0\to W_1\to\dots\to W_t$ such
that the words $W_0$ and $W_t$ are of length $n$, the maximal
length of $W_i$ is roughly $\exp n$, and $t$ is roughly $\exp\exp
n$. The corresponding diagram $\Delta$ for the conjugation of the
words $W_0$ and $W_t$, has area roughly equal to $ t\log t$. One
can obtain a diagram of area $t^2\log t/\log\log t$ and perimeter
$4t$ when gluing together $t/\log\log t$ copies of $\Delta$. This
proofs the property (2) of Theorem \ref{th1}.

To obtain the other inequalities, one has to strictly control what
the non-deterministic machine $M$ can do. (For example, the rules
of ages (2) and (5) look useless for property (2) but we need them
to prove the other properties.) The work of $M$ is studied in
Section 5. However one meets the primary difficulties when
proceeding to the calculation of areas for {\em arbitrary}
diagrams over the group $G$. In Section 6, we refine the technique
of \cite{OS}, and the exposition heavily depends on \cite{OS}. The
quadratic upper bound for the dispersion of a bipartite chord
diagram introduced in \cite{OS}, plays a key role here as well.

\section{Adding machine $Z(A)$}

Following \cite{OS}, we treat $S$-machines as HNN extensions of a
free group $F(Q,Y)$ generated by two sets of letters $Q=\cup_{i=1}^N
Q_i$ and $Y=\cup_{i=1}^{N-1} Y_i$ where $Q_i$ are disjoint and
non-empty (below we always assume that $Y_N=Y_0=\emptyset$). The set
$Q$ is called the set of $q$-letters, the set $Y$ is called the set
of $a$-letters.

Instead of the set of stable letters we have a collection $\Theta$
of $N$-tuples of  $\theta$-letters. Elements of $\Theta$ are
called {\em rules}. The components of $\theta$ are called {\em
brothers} $\theta_1,...,\theta_N$.

With every $\theta\in \Theta$, we associate two sequences of
elements in $F(Q\cup Y)$: $B(\theta)=[U_1,...,U_N]$,
$T(\theta)=[V_1,...,V_N]$, and a subsets $Y_i(\theta)\subseteq
Y_i$.

The generating set $\xxx$ of the group $S$ consists of all $q$-,
$a$- and $\theta$-letters. The relations, under condition
$\theta_{N+1}=\theta_1$,  are:

\begin{equation}\label{relations}
U_i\theta_{i+1}=\theta_i V_i,\,\,\,\, i=1,...,s, \qquad \theta_j
a=a\theta_j \;\;for\;\; all\;\; a\in Y_j(\theta)
\end{equation}

Sometimes we will denote the rule $\theta$ by $[U_1\to
V_1,...,U_N\to V_N]$. This notation contains all the necessary
information about the rule except for the sets $Y_i(\theta)$. In
most cases it will be clear what these sets are. By default
$Y_i(\theta)=Y_i$.

In this section we recall the definition and some properties of an
auxiliary adding machine $Z(A)$ introduced in \cite{OS}.

Let $A$ be a finite set of letters. Let the set $A_1$ be a copy of
$A$. It will be convenient to denote $A$ by $A_0$. For every
letter $a\in A$ $a_0$ and $a_1$ denotes its copy in $A_0$ and
$A_1$, respectfully.

The set of state letters of $Z(A)$ is $\{L\}\cup\{p(1), p(2),
p(3)\}\cup\{R\}$, i.e., there are 3 states for the $p$-letter, and
the letters $L$ and $R$ do not change their states. The set of
tape letters is $Y_1\cup Y_2$ where $Y_1=A_0\cup A_1$ and
$Y_2=A_0$.

The machine $Z(A)$ has the following rules (there $a$ is an
arbitrary letter from $A$) and their inverses. The comments
explain the meanings of these rules.

\begin{itemize}
\item $r_1(a)=[L\to L, p(1)\to a_1\iv p(1)a_0, R\to R]$.

\me {\em Comment.} The state letter $p(1)$ moves left searching
for a letter from $A_0$ and replacing letters from $A_1$ by their
copies in $A_0$.

\me

\item $r_{12}(a)=[L\to L, p(1)\to a_0\iv a_1p(2), R\to R]$.

\me

{\em Comment.} When the first letter $a_0$ of $A_0$ is found, it
is replaced by $a_1$, and $p$ turns into $p(2)$.

\me

\item $r_2(a)=[L\to L, p(2)\to a_0p(2)a_0\iv, R\to R]$.

\me

{\em Comment.} The state letter $p(2)$ moves toward $R$.

\me

\item $r_{21}=[L\to L, p(2)\tool p(1), R\to R]$, $Y_1(r_{21})=Y_1,
Y_2(r_{21})=\emptyset$.

\me

{\em Comment.} $p(2)$ and $R$ meet, the cycle starts again.

\me

\item $r_{13}=[L\tool L, p(1)\to p(3), R\to R]$, $Y_1(r_{13})=\emptyset,
Y_2(r_{13})=A_0$.

\me

{\em Comment.} If $p(1)$ never finds a letter from $A_0$, the
cycle ends, $p(1)$ turns into $p(3)$; $p$ and $L$ must stay next
to each other in order for this rule to be executable.

\item $r_{3}(a)=[L\to L, p(3)\to a_0p(3)a_0\iv, R\to R]$,
$Y_1(r_3(a))=Y_2(r_3(a))=A_0$

 \me

{\em Comment.} The letter $r_3$ returns to $R$.

\end{itemize}

\begin{rk} {If we replace every letter in $A_i$ by its index $i$,
then every word $u$ in the alphabet $A_0\cup A_1$ turns into a
binary number $b(u)$. If the machine starts with the word
$Lup(1)R$ where $u$ is a positive word in $A_0$, then $b(u)=0$ and
each 'regular' cycle of the machine adds $1$ to $b(u)$ modulo
$2^{|u|}$. After $2^{|u|}$ we obtain $Lup(3)R$.} \end{rk}

For every letter $a\in A$ we set $r_i(a\iv)=r_i(a)\iv$
($i=1,2,3$).

\begin{rk} {\rm All the rules of machine $Z(A)$ are transformed into
relations by formulas \ref{relations}, and so $Z(A)$ can also be
considered as a group. However, as in \cite{OS}, we will use
diagram and machine concepts in our proofs. All of them can be
found in \cite{OS}: reduced diagrams, bands in diagrams, trapezia,
their heights, bases and histories; admissible words, computation
$W=W_0\to_{\theta_1} W_1\to_{\theta_2} \dots\to_{\theta_t}
W_t=f\cdot W$ with history $f=\theta_1\theta_2\dots\theta_t$
determined by a trapezia, the width and area of a
computation.}\end{rk}

There is an obvious mirror analog $Z(A,mir)$ of the machine
$Z(A)$: $Y_1(mir)=A_0, Y_2(mir)=A_0\cup A_1$, and, for example,
the mirror analog of the rule $r_1(a)$ is $r_1(a,mir)=[L\to L,
p(1)\to a_0 p(1)a_1^{-1}, R\to R]$. (The state letter moves right
searching for a letter from $A_0$ and replacing letters from $A_1$
by their copies in $A_0$.) There are obvious mirror analogs of
lemmas \ref{lm901} - \ref{lm89}, but we will not formulate these
analogs. We often use $p$ instead of some $p(i)$ in subsequent
formulations.

\begin{lemma} \label{lm901} Suppose (1) $\base(W)\in\{LpR, p\iv pR\}$, both $W$ and
$f\cdot W$ contain $p(1)R$ (resp. $p(3)R$) or (2) $\base(W)$ is
$p\iv pR$ and both $W$ and $f\cdot W$ contain $p(1)R$ or $p(3)R$.
Assume that all $a$-letters in $W$ and in $f\cdot W$ are from
$A_0$ in both cases. Then $f$ is empty.
\end{lemma}

\proof In case (1), the assertion is proved in \cite{OS}, Lemma
2.27. It remains to consider case (2) where $W$ contains $p(1)$
but $f\cdot W$ contains $p(3)$. But this is impossible since the
$p$-letter cannot change state from $p(1)$ to $p(3)$ when its left
neighbor in the base differs from $L$. (See the definition of the
rule $r_{13}$.)
\endproof

\begin{lemma} \label{lm93} Let $W=LvpuR$, $\base(W)=LpR$. Suppose that
$|\theta\cdot W|>|W|$. Then for every computation $W\to_\theta
W_1\to W_2\to...\to f\cdot W$, we have $|W_i|>|W|$ for every $i\ge
1$.
\end{lemma}

\proof This assertion is proved in \cite{OS}, Lemma 2.24.
\endproof

\begin{lemma} \label{lm569} Suppose that an admissible word $W$ has the form
$LupvR$ (resp. $p\iv upvR$ or $Lvpup^{-1}$) where $u,v$ are words
in $(A_0\cup A_1)^{\pm 1}$. Let $\theta\cdot W=Lu'p'v'R$ (resp.
$\theta\cdot W=(p')\iv u'p'v'R$ or $\theta\cdot
W=Lv'p'u'p'^{-1}$). Then the projections of $uv$ and $u'v'$ (resp.
$v\iv uv$ and $(v')\iv u'v'$, or $vuv^{-1}$ and $v' u'(v')^{-1}$)
onto $A$ are freely equal.
\end{lemma}

\proof This assertion is proved in \cite{OS}, Lemma 2.18.
\endproof

\begin{lemma}\label{qqq} Suppose that one of the following conditions for an
admissible word $W$ of $Z(A)$ is satisfied (there $p=
\{p(1),p(2),p(3)\}$):
 $W$ does not contain a $p$-letter or
$\base(W)=Lpp\iv$, or $\base(W)=pp\iv p$, or $\base(W)=p\iv pR$,
or $\base(W)=LpR$. Then the width of any computation
$$W=W_0\to_{\theta_0} W_1\to_{\theta_1}...\to_{\theta_{t-1}} W_t$$
is at most $C\max(|W|,|W_t|)$ for some constant $C$.
\end{lemma}

\proof This assertion is proved in \cite{OS}, Lemma 2.29.
\endproof

\begin{lemma}\label{lm90} Let $\base(W)=LpR$.
Then for every computation $W=W_0\to W_1\to \dots\to W_t=f\cdot W$
of the $S$-machine $Z(A)$:

\begin{enumerate}
\item $|W_i|\le \max(|W|, |f\cdot W|)$, $i=0,...,t$,

\item If $W=LupR$ where $p=p(1)$ (resp. $p=p(3)$), $f\cdot W$
contains $p(3)R$ (resp. $p(1)R$) and all $a$-letters in $W, f\cdot
W$ are from $A_0^{\pm 1}$, then the length $g(|u|)$ of $f$ is
between $2^{|u|}$ and $6\cdot 2^{|u|}$, $u$ is a positive word,
and all words in the computation have the same length. Vice versa,
for every positive word $u$, such a computation does exist.
\end{enumerate}
\end{lemma}

\proof This assertion is proved in \cite{OS}, Lemma 2.25 and
Remark 2.19.
\endproof

\begin{lemma}\label{lm89} For every admissible word $W$ with $\base(W)=LpR$,
every rule $\theta$ applicable to $W$, and every natural number
$t>1$, there is at most one computation $W\to_\theta W_1\to ...\to
W_t$ of length $t$ where the lengths of the words are all the
same.
\end{lemma}

\proof This assertion is proved in \cite{OS}, Lemma 2.21.
\endproof

\section{How machine $\cal M $ works}

In this section, we introduce the machine ${\cal M}$ defining our
group $G$.

We set $N=5$ and ${\cal Q}=\cup_{i=1}^{5}{\cal Q}_i$. Here ${\cal
Q}_6={\cal Q}_1 =\{k_0\}$, ${\cal Q}_3 =\{k_1\}$, ${\cal Q}_5
=\{k_2\}$, ${\cal Q}_2$ is of the form $\{q_1^*(*)\}$, ${\cal
Q}_4$ is of the form $\{q_2^*(*)\}$, where the stars $^*$ and
$(*)$ will be replaced by particular indices below. The set of
rules of machine ${\cal M}$ will be partitioned in several {\em
ages}.

Consider the machine $Z(\{a\})$ for a 1-letter alphabet $\{a\}$.
Let $\Upsilon$ be the set of its rules. We introduce letters
$a_\tau\in A(\Upsilon)$ for every $\tau\in\Upsilon$.  Then we have
two copies $A(\Upsilon)_0$ and $A(\Upsilon)_1$ of this alphabet.
Let $Y = \cup_{i=1}^4 Y_i$ where $Y_1=\{a_0\}\cup \{a_1\}$, $Y_2
=\{a_0\}$, $Y_3= A(\Upsilon)_0$, $Y_4 = A(\Upsilon)_0\cup
A(\Upsilon)_1$.

{\bf Age}(1)  We correspond, to every rule $\tau$ of $Z(A)$, a
rule
  $\tau^1$ of {\em age } (1) of the machine ${\cal M}$. For example, for
  $\tau=r_1(a)$, we have
$$r_1(a)^1=[k_0\to k_0, q_1(1)^1\to a_1^{-1} q_1(1)^1 a_0, k_1\to k_1,
q_2^1\to a_{r_1(a)}q_2^1, k_2\to k_2]$$

$$ with \;\;Y_3(r_1(a)^1)=A(\Upsilon)_0, Y_4(r_1(a)^1)=\emptyset $$

\me {\em Comment.} Now the machine $Z(\{a\})$ works  with letters
$L, p, R$ replaced by $k_0, q_1, k_1$. At the same time it writes
the history of its work in alphabet $A(\Upsilon)_0$ on the tape
between the heads $k_1$ and $q_2$. For example, it can start
working with a word $k_0a_0^nq_1(1)^1k_1q_2^1k_2$ and finish
'adding' with $k_0a_0^nq_1(3)^1k_1uq_2^1k_2$ where $u$ is the
history $f$ of such a computation copied in the alphabet
$A(\Upsilon)$. The length of the positive word $u$ is $g(n)$ (see
Lemma \ref{lm90}).

{\bf Age}(12) The only {\em connecting} rule to the age (2) is

$$r^{12}=[k_0\to k_0, q_1(3)^1\to q_1(3)^2 , k_1\to k_1,
q_2^1\to q_2^2, k_2\to k_2]$$

$$ with\;\; Y_1(r^{12}) =\{a_0\}, Y_2(r^{12})=\emptyset,\;\;Y_3(r^{12})=A(\Upsilon)_0, Y_4(r^{12})=\emptyset $$

\me {\em Comment.} This rule changes the states of the heads $q_1,
q_2$ making possible the applications of rules of age (2). It is
applicable under the restrictions imposed on the sets
$Y_i(r^{12})$ above.

{\bf Age}(2) Again, we correspond to every rule $\tau$ of machine
$Z(\{a\})$ a rule $\tau^2$ of age (2) of machine ${\cal M}$. For
example, for $\tau = r_1(a)$, we define

$$r_1(a)^2=[k_0\to k_0, q_1(1)^2\to a_1^{-1} q_1(1)^2 a_0, k_1\to k_1,
q_2^2\to a_{r_1(a)}q_2^2a_{r_1(a)}^{-1}, k_2\to k_2]$$

$$ with \;\;Y_3(r_1(a)^2)=Y_4(r_1(a)^2)=A(\Upsilon)_0 $$

\me {\em Comment.} The work of ${\cal M}$ is similar to that in
age (1). But now the head $q_2$ runs to the left. For example, it
can start working with the word $k_0a_0^nq_1(3)^2k_1uq_2^2k_2$
(see Comment to Age (1)), then it can simulate the computation of
$Z(A)$ with history $f^{-1}$ and finish  'adding' with
$k_0a_0^nq_1(1)^2k_1q_2^2k_2$. We show and use that such a smooth
work between applications of rules of ages (12) and (23) is
possible only when the word $u$ has length $g(n)$ for some $n$.

{\bf Age}(23) The {\em connecting} rule to the age (3) is

$$r^{23}=[k_0\to k_0, q_1(1)^2\to q_1^3 , k_1\to k_1,
q_2^2\to q_2(1)^3, k_2\to k_2]$$

$$ with\;\; Y_1(r^{23}) =\{a_0\}, Y_2(r^{23})=\emptyset,\;, Y_3(r^{23})=\emptyset, Y_4(r^{23})=A(\Upsilon)_0,$$

\me {\em Comment.} The role of this rule is similar to that of
$r^{12}$.

{\bf Age}(3) Here we use the machine $Z(A(\Upsilon),mir)$. To
every rule $\tau$ of $Z(A(\Upsilon),mir)$, we correspond a rule
  $\tau^3$ of {\em age } (3) for the machine ${\cal M}$. For example, for
  $\tau=r_1(a)$ ($a\in A(\Upsilon)$) we have
$$r_1(a)^3=[k_0\to k_0, q_1^3\to q_1^3, k_1\to k_1,
q_2(1)^3\to a_0q_2(1)^3a_1^{-1}, k_2\to k_2]$$

$$ with \;\;Y_1(r_1(a)^3)=A(\Upsilon)_0, Y_2(r_1(a)^3)=\emptyset, Y_3(r_1(a)^3)=A(\Upsilon)_0 $$

\me {\em Comment.} The machine $Z(A(\Upsilon),mir)$ works now with
heads $L, p, R$ replaced by $k_1, q_2, k_2$. The head $q_1$ stays
by $k_2$, and the piece of tape between $k_1$ and $q_1$ is
unchanged. For example, it can start working with
$k_0a_0^nq_1^3k_1q_2(1)^3uk_2$ where $u$ is a reduced word of
length $g(n)$ in the alphabet $A(\Upsilon)$, and finish 'adding'
with $k_0a_0^nq_1^3k_1q_2(3)^3uk_2$ after application of $g(g(n))$
rules (double exponential in $n$ time by Lemma \ref{lm90}).

{\bf Age}(34) The connecting rule of age (34) is

$$r^{34}=[k_0\to
k_0, q_1^3\to q_1^4 , k_1\to k_1, q_2(3)^3\to q_2(3)^4, k_2\to
k_2]$$

$$ with\;\; Y_1(r^{34}) =\{a_0\}, Y_2(r^{34})=\emptyset,\;, Y_3(r^{34})=\emptyset, Y_4(r^{34})=A(\Upsilon)_0,$$

{\bf Ages} (4), (45), (5), (56), (6). The rules of ages (4), (5),
and (6) are similar to the rules of ages (3), (2) and (1),
respectively, up to the superscripts  at $r$- and $q$-letters: we
replace $1$ by $6$, $2$ by $5$, and $3$ by $4$. The connecting
rule of ages (45) and (56) are, respectively,

$$r^{45}=[k_0\to k_0, q_1^4\to q_1(1)^5 , k_1\to k_1,
q_2(1)^4\to q_2^5, k_2\to k_2]$$

$$ with\;\; Y_1(r^{45}) =\{a_0\}, Y_2(r^{45})=\emptyset,\;, Y_3(r^{45})=\emptyset, Y_4(r^{45})=A(\Upsilon)_0,$$
and

$$r^{56}=[k_0\to k_0, q_1(3)^5\to q_1(3)^6 , k_1\to k_1,
q_2^5\to q_2^6, k_2\to k_2]$$

$$ with\;\; Y_1(r^{56}) =\{a_0\}, Y_2(r^{56})=\emptyset,\;\;Y_3(r^{56})=A(\Upsilon)_0, Y_4(r^{56})=\emptyset $$

\me {\em Comment.} Let us start with the word $k_0a_0^n
q_1^3k_1q_2(3)^3 uk_2$. (See Comment to Age (3).) Then consecutive
applications of rules of ages (34), (4), (45), (5), (56) and (6)
can transform it as follows: $\to
k_0a_0^nq_1^4k_1q_2(3)^4uk_2\to\dots\to
k_0a_0^nq_1^4k_1q_2(1)^4uk_2\to
k_0a_0^nq_1(1)^5k_1q_2^5uk_2\to\dots\to
k_0a_0^nq_1(3)^5k_1uq_2^5k_2\to
k_0a_0^nq_1(3)^6k_1uq_2^6k_2\to\dots\to
k_0a_0^nq_1(1)^6k_1q_2^6k_2$. The computation
$k_0a_0^nq_1(1)^1k_2q_2^1k_2\to\dots\to
k_0a_0^nq_1(1)^6k_1q_2^6k_2$ we have considered as an example in
the comments to the definition of machine ${\cal M}$, has an
exponential width in $n$ and double exponential length of the
history.

As it was explained in the previous section, the machine ${\cal
M}$ defines the group $G=G({\cal M})$. (See (\ref{relations}).)

If the history $h$ of a computation is a product $h_1h_2\dots
h_s$, where, for every subword $h_i$, each of its letter has the
same age ($j_i$) ($j_i\in \{(1), (12),\dots, (56),(6)\}$), and
$j_i\ne j_{i+1}$ for $i=1,\dots,s-1$, then we say that this
computation has {\em brief history} $(j_1)(j_2)\dots (j_s)$.

Since a history $h$ is always a reduced word, it cannot contain
subwords of the form $\tau\tau^{-1}$. If $\tau$ is a connecting
rule, and the  computation base has at least one $q$-letter, the
$h$ has no subwords $\tau^2$, since every connecting rule changes
the states of $q$ letters. For the same reason, a brief history of
such a computation (with a $q$-letter in the base) cannot be of
the form $(1)(3)$ or $(6)(1)$, or $(3)(34)(3)$, etc.

We call a computation $W_0\too W_1\too\dots W_t$ {\em long} if the
brief history of this computation or of the inverse computation
has a subword equal to $(1)(12)(2)\dots(56)(6)$. Otherwise it is
{\em short}.

\section{$\cal M $-computations with various bases and histories}

 Denote by $\aaa$ the alphabet of all $a$-letters
  $\{a_0^{\pm 1},a_1^{\pm 1}\}\cup A(\Upsilon)_0^{\pm 1}\cup A(\Upsilon)_1^{\pm 1}$   .
The length of a word $W$ we denote by $|W|$, and the $a$-{\em
width} $|W|_a$ of $W$ is the number of $a$-letters in $W$. The
width $||W||$ of an admissible word of the form $q^{-1}uqvk$ and
$kvquq^{-1}$ where $u$ and $v$ are words in $\aaa$, is defined as
$3+|u|+2|v|$, and $||W||=|W|$ for all other words by definition.

We say that a reduced computation is {\em regular} if the
applications of its rules do not change the width. An application
of a rule $W\too W'$ {\em increases the width} of $W$ if the word
$W'$ is longer than $W$.  Then the application of the inverse rule
to $W'$ {\em decreases the width}.

\begin{lemma}\label{decrincr} Let $W_0\too_{\rho_1}\dots\too_{\rho_l} W_l$ be a
computation with base $k_1q_2k_2$. Assume that all the rules
$\rho_1,\dots,\rho_l$ are of age (1) or (6) ( of age (2), or (5)).
Then there is an integer $d$, $0\le d \le l$,  such that the
applications of $\rho_1,\dots,\rho_d$ decrease (do not increase)
the widths of words $W_0,\dots,W_{d-1}$, and the applications of
$\rho_{d+1}, \dots,\rho_l$ increase (do not decrease) the widths
of $W_d,\dots,W_{l-1}$.
\end{lemma}

\proof Assume that an application of a rule $\rho_i$ of age (1)
increases the width of $W_{i-1}=k_1wq_2k_2$, and $i<l$. Then
$W_i=k_1waq_2k_2$ with a reduced word $wa$, $a\in
A(\Upsilon)_0^{\pm 1}$. Since the history of a computation is
reduced, we have $\rho_{i+1}\ne \rho_i^{-1}$, and so the
application of $\rho_{i+1}$ must also increase the width of $W_i$
as this follows from the definition of the rules of age (1). The
lemma statement follows from this observation. The proofs of the
assertion for rules of ages (6), (2) and (5), are similar.
\endproof

\begin{lemma}\label{vozvrat} Let the history of a computation
be $\eta h \eta^{-1}$, where $\eta$ is a connecting rule and $h$
has no connecting rules. Assume that the base of this computation
has one of the forms $kqq^{-1}$, $qq^{-1}q$, $kqk$, $q^{-1}qk$.
Then the base is $q^{-1}qk$ or $kqk$, and if $q=q_1$, then rules
of $h$ have age (3) or (4), and if $q=q_2$, then rules of $h$ have
age (1) or (6).
\end{lemma}

\proof No  connecting rule is applicable to a word $quq^{-1}vq$
where $u$ and $v$ are $a$-words.

Now we assume that the base is $kqk$ and the age of $h$ is not (3)
or (4). Then the equality $q=q_1$ is impossible by Lemma
\ref{lm901}. Similarly, $h$ cannot be of age (3) or (4) if
$q=q_2$. The assumption that the history $h$ is of age (2) or (5)
and the base is $k_1q_2k_2$ leads to a contradiction since $h$ is
reduced, and a connecting rule $\rho$ is applicable when
$Y_4(\rho)=\emptyset$. Lemma \ref{lm901} also works for bases
$kqq^{-1}$, $q^{-1}qk$ if a connecting rule is applicable to a
word having such a base.
\endproof

\begin{lemma}\label{2conrules} Let the base of a computation $W_0\too
W_1\too\dots W_t$ have one of the forms $kqq^{-1}$, $k^{-1}k$,
$qq^{-1}q$, $kqk$, $q^{-1}qk$, $kk^{-1}$, or $k_2k_0$. Then

(1) all applications of the rules are regular if the first and the
last rules are both connecting rules and there are no subwords
(12)(1)(12) and (56)(6)(56) in the brief history; there are no
such subwords if the base contains $q=q_1$;

(2) if there is an application $W_{i-1}\too W_i$ of a connecting
rule in the computation and there are no letter $q_2$ in the base
or there are no rules of ages (1) and (6) in the history, then
$||W_i||\le ||W_s||$ for arbitrary $s\in\{0,1,\dots,t\}$.
\end{lemma}

\proof (1) We may exclude cases $kk^{-1}$, $k^{-1}k$ and $k_2k_0$
since they are trivial: $||W_0||=||W_1||=\dots$. Then we may
assume that the only connecting rules are the first one and the
last one. The remaining rules must be of the same age, say $(l)$
where $l\in \{2,3,4,5\}$. (If, for example, $l=1$, then there must
be a subword $(12)(1)(12)$ in the brief history, and $q=q_2$ by
Lemma \ref{vozvrat}.) The $a$-letters of both $W_0$ and $W_t$ must
belong to $\{a_0^{\pm1}\}$ or to $A(\Upsilon)_0^{\pm 1}$.

Assume the base is $kqk$.  It is easy to see that the applications
of rules  do not change the projection of a word onto the
subalphabet $A_0$ since neither of the rules are of age (1) or
(6). It follows that $||W_0||=||W_t||$ and $||W_i||\ge ||W_0||$
for $i=1,\dots,t-1$. We notice now that no rule application is
increasing by Lemma \ref{lm93}, if the age $(l)$ is (2) or (5) and
the base is $k_0q_2k_1$, or $l=3,4$ and the base is $k_1q_1k_2$
(i.e., a copy of machine $Z(\{a\})$ or machine $Z(A(\Upsilon),
mir)$ works). In other cases, the assertion follows from Lemma
\ref{decrincr} since $||W_0||=||W_t||$.

Let the base be $kqq^{-1}$. Then the connecting rules (12) and
(56) are not applicable, and one may assume by the symmetry that
$l=3$. The mirror version of Lemma \ref{lm901} makes this case
impossible if $q=q_2$. Otherwise we just have $W_0=W_1=\dots$.

Similar arguments work for the bases $q^{-1}qk$ and $qq^{-1}q$.

(2) Let $W_s=k_0u_sq_1v_sk_1$. The
  reduced form of the projections of $u_sv_s$ on the
  alphabet $A_0$ do not depend on $s$ by Lemma \ref{lm569}. But
  $u_iv_i$ is a reduced word in $A_0$ since one of the factors is
  empty (recall that $W_i$ is the result of an application of a
  connecting rule). Hence $||W_s||\ge ||W_i||$. The argument is
  similar if $q=q_2$ and there are no rules of ages (1) and (6) in
  the history.

Again, the base cannot be equal to $qq^{-1}q$, and the statement
is obvious for bases $kk^{-1}$, $k^{-1}k$ and $k_2k_0$. Then we
may assume by part (1), that $i=1$, the history is $\eta h$ where
$\eta$ is a connecting rule and $h$ has no connecting rules.

  Assume, for example, that the base is $q_1^{-1}q_1k_1$, the rule
  $\eta$ is of age (12) or (23) ((45),  or (56)) and $h$ is of age (2) (of age (5)).
  Let $W_s=q_1^{-1}u_sq_1v_sk_1$. Then the
  reduced forms of the projections of $v_s^{-1}u_sv_s$ on the
  alphabet $A_0$ do not depend on $s$ by Lemma \ref{lm569}. Since
  $u_1$ is a word in $A_0$ and $v_1$ is empty (recall that $\eta$ is a
  connecting rule of age (12) or (23)), we have
  $$||W_s||=3+|u_s|+2|v_s|\ge 3+|v_s^{-1}u_sv_s|\ge 3+|u_1|=||W_1||$$
  as desired.  If $h$ is of age (3) or (4), then obviously
  $||W_0||=||W_1||=\dots$.

  Similar arguments work for $q=q_2$ and also for bases of the
  form $kqq^{-1}$. The lemma is proved.
\endproof

\begin{lemma} \label{short} Let the base of a computation $W_0\too W_1\too\dots\too W_t$
be one of the forms $kqq^{-1}$, $k^{-1}k$, $qq^{-1}q$, $kqk$,
$q^{-1}qk$, $kk^{-1}$. Then

(1) if the history of the computation contains both connecting
rules $r^{12}$ and $r^{23}$ (or their inverses), then the base has
the form $kqk$;

(2) if $q=q_1$ or the computation is short, then $|W_i|\le c
\max(|W_0|,|W_t|)$ for some constant $c$ independent of the
computation.
\end{lemma}

\proof (1) The connecting rule $r^{12}$ is not applicable whenever
$kqq^{-1}$ is the base.

We have $q=q_1$ if the base is $q^{-1}qk$, since otherwise
$r^{23}$ is not applicable. But after the application of $r^{12}$,
the state of the $q$-letters is $q_1(3)^2$, and before the rule
$r^{23}$ is applied, the state must be $q_1(1)^2$. But the state
$q_1(1)^2$ cannot be reached since the base has a subword
$q_1^{-1}q_1$ (but not $k_0q_1$), a contradiction.

Similarly, the bases of forms $k^{-1}k$, $qq^{-1}q$, $kqq^{-1}$,
and $kk^{-1}$ can be eliminated.

(2) We can assume that there is a $q$-letter in the base since
otherwise the assertion is obvious. Let $h=\rho_0\dots\rho_t$ be
the computation history, and $\rho_{i_1}, \dots , \rho_{i_l}$ all
connecting rules in this computation.

For $l=0$,  the assertion follows from Lemma \ref{qqq} if the age
and the base are similar to those for machines $Z(\{a\})$ or
$Z(\Upsilon), mir)$. Otherwise it either obvious or follows from
Lemma \ref{decrincr}.

 Let $l\ge 1$. Then denote by $h_0,h_1,\dots, h_l$ the subwords of $h$ such that
$h_j$ starts with $\rho_{i_j}$ (with $\rho_0$ for $j=0$) and
terminates with $\rho_{i_{j+1}}$ (with $\rho_t$ for $j=l$).

First assume  that either $q=q_1$ or $h_0$ has no rules of age (1)
or (6). Then by Lemma \ref{2conrules}(1), we have that
$||W_{i_1-1}||=||W_{i_l}||=||W_s||$ for $i_1-1\le s\le i_l$, and
by Lemma \ref{2conrules}(2), $||W_0||\ge||W_{i_1}||=||W_{i_l}||$.
Therefore this case is reduced to the statement for the
subcomputations $W_0\too\dots\too W_{i_1-1}$ and
$W_{i_s}\too\dots\too W_t$ having no connecting rules in the
histories. The case where $h_l$ has no rules of age (1) or (6) is
similar.

Thus, it remains to eliminate the case: $l\ge 1$,  $q=q_2$, $h_0$
contains rules of age (1) or (6), and similarly $h_l$ does.
Therefore $l>1$, and we may assume that $h_0$ contains a rule of
age (1). Then $\rho_{i_1}=(12)$, and it follows from Lemma
\ref{vozvrat} that $\rho_{i_2}=(23)$. Then the base is $kq_2k$ by
part (1) of the lemma. Now applying Lemma \ref{vozvrat} several
times, we have that $\rho_{i_3}=(34)$, $\rho_{i_4}=(45)$,
$\rho_{i_5}=(56)$, $l=5$, and the computation is long against the
lemma condition.

\endproof

\begin{lemma}\label{log} Let $W_0\too W_1\too\dots\too W_t$ be
a long computation with base $k_0q_1k_1q_2k_2$ or
$k_1q_2k_2k_0q_1k_1$, or $k_2k_0q_1k_1q_2k_2$. Then for every
$W_i$, where $1\le i\le t$, we have $|W_i|_a\le\max(|W_1|_a,
|W_t|_a, 2log_2 t)$. If $k_0uq_1$ is a subword of $W_0$ with
$|u|=n$, and a subcomputation $W_l\too W_{l+1}\too\dots\too W_m$
starts (ends) with an application of the rule (12) (the rule
(56)), then $m-l = 5+2g(n)+2g(g(n))$ for some integer $n$, and the
$a$-width of this subcomputation is $n+g(n)$. The $a$-widths of
the restrictions of this subcomputation to bases $k_0q_1k_2$ and
$k_1q_2k_2$ are $n$ and $g(n)$, respectively.
\end{lemma}

\proof We will assume that the base is $k_0q_1k_1q_2k_2$. By Lemma
\ref{2conrules} (1), there are no subwords (12)(1)(12) and
(56)(6)(56) in the brief history $B$ of the computation, and so
$B=(1)(12)\dots (56)(6)$.

Denote by $W_i'$ and $W''_i$, respectively, the prefix (the
suffix) of $W_i$ ending (starting) with $k_1$. Then $|W'_i|\le
\max (|W'_0|, |W'_l|)$ for $0\le i\le l$, by Lemma \ref{lm90}, and
$|W'_l|\le |W'_0|$ by Lemma \ref{2conrules} (2). Also we have
$|W''_i|\le \max (|W_0''|, |W_l''|)$ by Lemma \ref{decrincr}.

Then $|W_i| = |W_l|=|W_m|$ for $l\le i\le m$ by Lemma
\ref{2conrules} (1), and $t\ge g(|v|_a)\ge 2^{|v|_a}$ by Lemma
\ref{lm90}  for the subcomputation of age (3), where $v$ is any of
the words read between $k_1$ and $k_2$ in age (3). Similarly, when
considering the maximal subcomputation of age (2), we have by
lemmas \ref{2conrules} and \ref{lm90}, $|v|_a =g(|u|)=g(n)\ge
2^n$. Hence, for $l\le i\le m$, we have $|W_i|_a\le n + |v|<
\log_2\log_2 t+\log_2 t\le 2\log_2 t$.  For $i<l$, $|W_i|_a=
|W_i'|_a+|W_i''|_a \le \max (|W_0|_a, n +|v|_a) \le \max (|W_0|_a,
2log_2t)$. Since there is a similar estimate in case $i\ge m$, the
desired upper bound for $|W_i|_a$ is obtained for all $i$.

The equality $m-l = 5+2g(n)+2g(g(n))$ follows from lemmas
\ref{lm89} and \ref{lm90} since there are 5 connecting rules in
the subcomputation.
\endproof

\begin{lemma}\label{area-kqkqk} Let $W_0\too W_1\too\dots\too W_t$ be
a computation with base $k_0q_1k_1q_2k_2$ or $k_1q_2k_2k_0q_1k_1$,
or $k_2k_0q_1k_1q_2k_2$. Assume that $10g(g(n-1))\le t\le g(g(n))$
for some integer $n$. Then the area of corresponding trapezium
$\Delta$ does not exceed $Ct(|W_0|_a+|W_t|_a)$ for a constant $C$
independent of the computation.
\end{lemma}

\proof If the computation is short, then the statement follows
from Lemma \ref{short} applied to the restrictions of the
computation to subbases $k_0q_1k_1$ and $k_1q_2k_2$. Therefore as
in the proof of Lemma \ref{log}, one may suppose that the brief
history of the computation is $(1)(12)\dots(6)$.

We denote by  $T_i$ the $i$-th $\theta$-band of $\Delta$. Observe
that at most four $(\theta,a)$-cells of $T_i$ can be attached to
its $(\theta, q)$-cells along $a$-edges. Hence the number of cells
in $T_i$ does not exceed $4+6+|W_i|=10+|W_i|$ because $T$ has at
most six $k$- and $q$-cells.

Let the application of rules (12) and (56) be the $l$-th and the
$m$-th, respectively, in the history $\rho_1\dots\rho_t$. Then
$\Delta$ is the union of of 3 subtrapezia $\Delta_1$, $\Delta_2$
and $\Delta_3$ which correspond to subwords $\rho_1\dots\rho_{l}$,
$\rho_l\dots\rho_m$, and $\rho_{m}\dots\rho_t$, respectively. By
Lemma \ref{log}, there is an integer $r$ such that
$m-l=5+2g(r)+2g(g(r))$, and $r<n$ since $m-l\le t<g(g(n))$. Also,
by Lemma \ref{log}, $|W_i|_a\le r+g(r)$ for $l\le i\le m$. Using
the observation of the previous paragraph, we see that the area of
$\Delta_2$ does not exceed $(10+r+g(r))(5+2g(r)+2g(g(r)))$.

Consider the restriction of the subcomputation with subhistory
$\rho_1\dots\rho_l$ to the subbase $k_1q_2k_2$. By Lemma
\ref{decrincr}, there is an integer $d$, $0\le d \le l$, such that
the applications of $\rho_1,\dots,\rho_d$ decrease the $a$-widths
of the suffices $V_{i-1}=k_1\dots k_2$ of subwords $W_{i-1}$, and
the applications of $\rho_{d+1}, \dots,\rho_l$ increase them. But
$|V_l|_a=g(r)$ by Lemma \ref{log}, and therefore $l-d\le g(r)$.

Let $\Delta_{11}$ and $\Delta_{12}$ be the subtrapesia of
$\Delta_1$ of heights $d$ and $l-d$, respectively. It follows from
Lemma \ref{log} that the $a$-width of $\Delta_{12}$ does not
exceed $r+g(r)$, and therefore its area does not exceed
$(10+r+g(r))g(r)$ since its height is not greater than $g(r)$.

The $a$-width of $\Delta_{11}$ is not greater than $3l+r+g(r)$
because a single application of a rule changes the $a$-width of a
word with base $k_0q_1k_1q_2k_2$ at most by $3$. Thus the area of
$\Delta_{11}$ does not exceed $l(10+3l+r+g(r)).$ Therefore the
area of $\Delta_1$ is not greater than
$(10+r+g(r))g(r)+l(10+3l+r+g(r))$. Similarly, the area of
$\Delta_3$ is bounded from above by
$(10+r+g(r))g(r)+(t-m)(10+3(t-m)+r+g(r))$. Thus the area of
$\Delta$ is at most $100g(g(r))g(r)+100tg(r)+3t^2.$

It follows from the definition of $\Delta_{11}$ that $|W|_a \ge
d$. Also recall that $d\ge l-g(r)$. Hence $|W_0|_a+|W_t|_a$ is at
least
$$l-g(r) + (t-m) -g(r) = t -(m-l) - 2g(r)= t-(5+2g(r)+2g(g(r)))-2g(r)>\max (t/3,
g(g(r)))$$
 by the choice of $t$ and by inequality $r\le n-1$. Therefore the area of
$\Delta$ is not greater than $Ct(|W_0|_a+|W_t|_a)$ for a constant
$C$ independent of the computation.
\endproof

\begin{lemma}\label{loglog} Let $W_0\too W_1\too\dots\too W_t$ be
a long computation with base $k_0q_1k_1q_2k_2$ or
$k_1q_2k_2k_0q_1k_1$, or $k_2k_0q_1k_1q_2k_2$. Then
$|W_0|_a+|W_t|_a+1\ge c_0 \log'\log't$ for a positive constant
$c_0$.
\end{lemma}

\proof We will use the notation of Lemma \ref{area-kqkqk}. To
proof the statement, it suffices to assume that $d=0$ because, by
Lemma \ref{lm90} (2), the applications of $\rho_1,\dots,\rho_d$
cannot increase the lengths of subwords $k_0\dots k_1$ of the
words $W_i$ since the rules of age (1) follows by a connecting
rule in the whole computation; and they  decrease the lengths of
their subwords of the form $k_1\dots k_2$. Similar assumption is
applicable to $W_m\too W_{m+1}\too\dots\too W_t$. Then $l=l-d\le
g(r)$ and $t-m\le g(r)$ as in the proof of Lemma \ref{area-kqkqk}.
Hence $t\le 5+2g(r)+2g(g(r))+2g(r)$.

On the other hand, by Lemma \ref{2conrules} (2) restricted to the
base $k_0q_1k_1$, we have $|W_0|_a, |W_t|_a\ge r$. To finalize the
proof, it suffices to note that, by Lemma \ref{lm90}, there exists
a positive $c_0$ such that $c_0\log'\log'(5+4g(r)+2g(g(r)))\le
2r+1$ for every $r\ge 0$.
\endproof

\begin{lemma}\label{width} Let $W_0\too W_1\too\dots\too W_t$ be
a long computation with  base $k\dots k$, where the first and the
last $k$-letters coincide. Then, for some constant $c$ and for
every $i\in\{1,\dots,t\}$, we have $|W_i|_a\le
c(|W_1|_a+|W_t|_a+b\log_2 t)$, where $b$ is the length of the
base.
\end{lemma}

\proof Neither the base nor its inverse word has subwords of the
form $kqq^{-1}$ or $k^{-1}k$, or $qq^{-1}q$, or $q^{-1}qk$, or
$kk^{-1}$ by Lemma \ref{short}. Recall also that every letter
$k_2$ can be followed in the base by letter $k_0$ only. Then it
follows from the lemma assumption that every word $W_i$ (or the
inverse word) can be covered by its subwords with base of the form
$k_0q_1k_1q_2k_2$ or $k_1q_2k_2k_0q_1k_1$, or
$k_2k_0q_1k_1q_2k_2$, in such way that every basic letter is
covered at most two times and every $a$-letter is covered once.
Now the assertion is a consequence of Lemma \ref{log}.
\endproof

\begin{lemma}\label{area} Let $\Delta$ be a trapezium of height $h\ge 1$ with
either (a) base $k\dots k$, where the first and the last
$k$-letters coincide, and which contains neither subwords
$(qq^{-1}q)^{\pm 1}$ nor shorter subwords of the form $k\dots k$,
or (b) base $qq^{-1}q$. Then the area of $\Delta$ does not exceed
$ch(|W|_a+|W'|_a+ \log' h)$ for a constant $c$, where $W, W'$ are
the labels of its top and bottom, respectively. The third summand
can be replaced by 1 if the base is $qq^{-1}q$ or $\Delta$
corresponds to a short computation.
\end{lemma}

\proof (a) It follows from the lemma assumption that the length
$b$ of the base is bounded from above. Since the area of the
$i$-th band of $\Delta$ can exceed the length of $W_i$ at most by
$2b$, the lemma statement follows from lemmas \ref{width} and
\ref{short}.

(b) Similarly, the statement follows from Lemma \ref{short} in
this case.
 \endproof

\section{Areas of diagrams over the group $G$}

As in \cite{OS}, we use constants $L,K,\delta$. I suffices to set
$L\ge6$ since the are no defining relations of length $>6$ now.
Then $K=2K_0$, where $K_0$ bounds from above the length of bases
having neither subwords $(qq^{-1}q)^{\pm 1}$ nor subwords $xux$
whith a $k$-letter $x$. As in \cite{OS}, a sufficiently small
positive $\delta$ is selected so that $\delta(4L+LK+1)<2$. As in
Section 4.1 \cite{OS}, the {\em lengths} of words, paths and
perimeters of diagrams are modified now. (The number of edges in a
path is called now a {\em combinatorial length}.) The reader of
this section should have the paper \cite{OS} at hand. In
particular, the concept of diagram dispersion ${\cal E}(\Delta)$
is crucial for the proofs of Lemma 6.2 \cite{OS} and the lemmas of
this section. However it is not defined here since we do not use
the definition and use the same property of dispersion as in
\cite{OS} (e.g., the quadratic upper bound in term of the
perimeter $|\partial\Delta|$). As in \cite{OS} we take a big
enough constant $M$. Here ``big enough" means that $M$ satisfies
the inequalities used in the proof of lemmas \ref{main2} and
\ref{main1}. Each of them has the form $M>C$ for some constant $C$
that does not depend on $M$ (but depends on the constants
introduced earlier). Since the number of inequalities is finite,
one can choose such a number $M$.

\begin{lemma}\label{main2} The area of a reduced diagram $\Delta$
does not exceed $M\Psi(n) +M\psi(n)\eee(\Delta)$, where
$n=|\partial\Delta|$ and $\psi(n)=log'n/log'log'n$.
\end{lemma}

\proof We follow the proof of Lemma 6.2 \cite{OS}. Steps 1 and 2
are analogous to those in \cite{OS}: The only difference is that
one must multiply the entropy ${\cal E}(\Delta)$ by $\psi(\Delta)$
and replace the factor $\log'(\dots)$ by $\psi(\dots)$. Then we
use all the notations of Step 3 \cite{OS} for the supposed minimal
counter example: $\Delta$, ${\cal T}$, ${\cal T}'$, ${\cal Q}$,
${\cal Q}'$, ${\cal Q}_2 - {\cal Q}_4$, $l$, $l'$,($l'>l/2$),
$l_3$, $l_4$, $\Gamma$, $\Gamma'$, $\Gamma_1 - \Gamma_4$,
$\Delta_0$, $n$, $n_0$, $\alpha_i$, $p_i$, $p^i$, $u_i$ $d_i$,
$d'_i$ for $i=3,4$ and $A_0 - A_4$. Then reader can just compare
our arguments here and there. In particular, as in (6.23)
\cite{OS},

\begin{equation}\label{nbezn0}
n-n_0\ge2+\delta(\max(0,d'_3-2L,
\alpha_3-(d_3-d'_3)-2Ll_3)+\max(0,
d'_4-L,\alpha_4-(d_4-d'_4)-2Ll_4))
\end{equation}

 By Lemma \ref{area} we have now
\begin{equation}\label{gamma22}
A_2\le C_2l'(d_3+d_4+\log' l')
\end{equation}
for some constant $C_2$, and

\begin{equation}\label{gamma22'}
A_2\le C_2l'(d_3+d_4+1)
\end{equation}
if the computation defined by the trapezium $\Gamma_2$ is short.

The inequalities (\ref{gamma22}), (\ref{gamma22'}) provide us with
the following modification of our task (in comparison with (6.28)
in \cite{OS}): To obtain the desired contradiction, we must now
prove that

\begin{equation}\label{tsel'2}
(Mn(n-n_0)+\frac{M}{K^2} l'(l-l'))\psi(n)\ge  C_3l'(d_3+d_4+\log'
l')+C_3(l_3^2+l_4^2)+2\alpha_3 l_3 + 2\alpha_4l_4
\end{equation}
where (as in inequality (6.28), \cite{OS}) $C_3\ge C_2$ is a
constant that does not depend on $M$, and if the trapezium
$\Gamma_2$ corresponds to a short computation, we must prove
(\ref{tsel'2}) with the logarithmic summand replaced by $1$.

First, as in \cite{OS}, we can choose $M$ big enough so that
\begin{equation}\label{MC3l2}
\frac{M}{3K^2}l'(l-l')\ge C_3(l_3^2+l_4^2)
\end{equation}

Then, as in \cite{OS}, we assume without loss of generality that
$\alpha_3\ge\alpha_4$, and consider two cases.

(a) Suppose we have $\alpha_3\le 2C_3(l-l')$.

Since $d_i\le \alpha_i+d'_i$ for $i=3,4$, we also, by inequality
(\ref{nbezn0}), have $d_3+d_4+1\le\alpha_3+\alpha_4+d'_3+d'_4+1<
4C_3(l-l')+\delta^{-1}(n-n_0)+2L-2\delta^{-1}+1 <
4C_3(l-l')+\delta^{-1}(n-n_0)$ because $\delta^{-1}>L+1/2$ by the
choice of $\delta$. Therefore

\begin{equation}\label{a122}
(\frac{M}{5K^2} l'(l-l')+\frac{M}{2}  n(n-n_0))\ge
C_3l'(d_3+d_4+1)
\end{equation}
because $n\ge l'$, $n-n_0\ge 2$ by (\ref{nbezn0}), $M\ge
C_3\delta^{-1}$ and $M\ge 20 C^2_3K^2$.

Since $l_3+l_4=l-l'<l'$, we have also
\begin{equation}\label{a322}
\frac{M}{5K^2} l'(l-l')\ge 2(l-l')2C_3(l-l')\ge
2\alpha_3l_3+2\alpha_4l_4
\end{equation}
because $M\ge 20K^2C_3$.

If the trapezium $\Gamma_2$ corresponds to a short computation,
then the inequality (\ref{tsel'2}) (with the logarithmic summand
replaced by $1$) follows from (\ref{MC3l2}), (\ref{a122}) and
(\ref{a322}). Then we assume that the computation is long. Since
the base of $\Gamma_2$ satisfies the Lemma \ref{area} condition
(as in \cite{OS}), it follows from Lemma \ref{short} that the base
or its inverse has one of the forms $k_0q_1k_1q_2k_2k_0$,
$k_1q_2k_2k_0q_1k_1$, $k_2k_0q_1k_1q_2k_2$.

Now we consider two possibilities.

(a1) Let $l-l'\ge \frac1{20C_3} \log'\log'l'$. Then
\begin{equation}\label{a42}
\frac{M}{5K^2} l'(l-l')\psi(n)\ge C_3l'\log'l'
\end{equation}
by the definition of the function $\psi(n)$, because $M\ge
100K^2C_3^2$. By adding inequalities (\ref{MC3l2}), (\ref{a122}) -
(\ref{a42}), we obtain a stronger inequality than the desired
inequality (\ref{tsel'2}).

(a2) Let $l-l'\le\frac1{20C_3}\log'\log'l'.$  Let us estimate the
number of $\aaa$-edges lying on the path $u_3$. Recall that it is
equal $d'_3$ plus $|p^3|_a$ (see \cite{OS}). It follows from Lemma
4.10 \cite{OS} that $|p^3|_a\ge (d_3-d'_3)-C_0l_3$ for a constant
$C_0$ (One may assume that $C_3>C_0c_0^{-1}/10$ where $c_o$ is
given by Lemma \ref{loglog}.) Thus $|u_3|_a\ge d_3 -C_0l_3$. By
using a similar lower bound for $|u_4|_a$, we have
$|u_3|_a+|u_4|_a\ge (d_3+d_4)-C_0(l-l')$. When applying the
assumption (a2) and Lemma \ref{loglog} to the right-hand side of
this equality, we have $|u_3|_a+|u_4|_a+1\ge
(c_0-\frac{C_0}{20C_3})\log'\log'l'$. By Lemma 4.6 \cite{OS}, we
obtain $|u_i|\ge l_i+\delta(|u_i|_a-Ll_i)$ for $i=3,4$. Since we
may chose $C_3$ such that $C_3\ge c_0^{-1}L$, and
$l_3+l_4=l-l'\le\frac{1}{20C_3}\log'\log' l'$, we have now
$$|u_3|+|u_4|+1\ge
l_3+l_4+\delta(c_0-\frac{C_0}{20C_3}-\frac
L{20C_3})\log'\log'l'\ge(l-l')+\frac{\delta c_0 }{2}\log'log'l'$$
It follows from this inequality and the comparison of perimeters
$|\Delta|$ and $|\Delta_0|$, that $n-n_0\ge
2+|u_3|+|u_4|-(l-l')\ge\frac{\delta c_0}{2}\log'log'l'$, and since
$l'\le n$ and $M\ge 2\delta^{-1}c_0^{-1}C_3$, we obtain

\begin{equation}\label{a52}
\frac{M}2 n(n-n_0)\psi(n)\ge C_3 l' log'l'
\end{equation}

 The sum of inequalities (\ref{MC3l2}), (\ref{a122}), (\ref{a322}),  and (\ref{a52}) gives us
 a stronger inequality than (\ref{tsel'2}).

(b) Assume now that $\alpha_3>2C_3(l-l')$. Then, as in \cite{OS},
we have
\begin{equation}\label{d34}
d_3+d_4+1\le \frac53\alpha_3+\delta^{-1}(n-n_0)
\end{equation}
Here we add $1$ to the left-hand side (comparatively to
\cite{OS}). This is possible since, as in \cite{OS}, $n-n_0\ge 2$
and $\delta$ can be selected small enough. Then, as in \cite{OS},
one ontains
\begin{equation}\label{raznitsa}
n-n_0\ge \frac17 \delta\alpha_3.
\end{equation}

From inequalities (\ref{d34}),(\ref{raznitsa}), $M\ge
10C_3\delta^{-1}$ and $l'\le n/2$, we have
\begin{equation}\label{b122}
\frac{M}3 n(n-n_0) \ge \frac{10C_3}3 \delta^{-1}n(n-n_0)\ge
C_3l'(d_3+d_4+1)
\end{equation}
Inequalities (\ref{raznitsa}), $M\ge 21\delta^{-1}$,
$\alpha_3\le\alpha_4$, and $l_3+l_4=l-l'\le \l/2\le\frac14 n$ give
us
\begin{equation}\label{b222}
\frac{M}6 n(n-n_0)\ge \frac{7}2 \delta^{-1}(n-n_0)n\ge
2\alpha_3(l_3+l_4)\ge 2\alpha_3l_3+2\alpha_4 l_4
\end{equation}

If the $\Gamma_2$-computation is short then the corresponding
version of (\ref{tsel'2}) (where the logarithmic summand is
replaced by $1$) follows from inequalities (\ref{MC3l2}),
(\ref{b122}) and (\ref{b222}). Thus, as in case (a), by Lemma
\ref{short}, the base or its inverse can be supposed having one of
the forms $k_0q_1k_1q_2k_2k_0$, $k_1q_2k_2k_0q_1k_1$,
$k_2k_0q_1k_1q_2k_2$.

 Then, from (\ref{d34}), (\ref{raznitsa}) and Lemma \ref{loglog}, we have
$$n-n_0\ge\frac{\delta}{13}(d_3+d_4+1)\ge\frac{\delta
c_0 }{13}\log'\log'l'$$ Since $M\ge 13\delta^{-1}c_0^{-1}C_3$ and
$l'\le n/2$, it follows from the definition of $\psi(n)$ that
\begin{equation}\label{b322}
\frac{M}2 n(n-n_0)\psi(n)\ge C_3 l' log'l'
\end{equation}

The inequality (\ref{tsel'2}) follows now from inequalities
(\ref{MC3l2}), (\ref{b122}), (\ref{b222}),  and (\ref{b322}).

 The lemma is proved by contradiction.
\endproof

\begin{lemma}\label{main1} Let the perimeter $n$ of a reduced diagram $\Delta$
satisfy inequality $ n \le g(g(r))$ for some positive integer $r$.
Then the area of diagram $\Delta$ does not exceed $M(n^2 +
m^2log'm) +M\eee(\Delta)$, where $m=10g(g(r-1))log'n$.
\end{lemma}

\proof By Lemma \ref{main2}, the statement is true for $n\le
10g(g(r-1))log'n$ since $m\ge n$ in this case. Then arguing by
contradiction, we consider a counter-example $\Delta$ with minimal
perimeter $n>10g(g(r-1))log'n$.

  As in the proof of Lemma \ref{main2}, we follow the proof of Lemma
6.2 \cite{OS}. Steps 1 and 2 are analogous to those in \cite{OS}
since the extra term $m^2\log'm$ does not affect. Then again we
use the notations of Step 3 \cite{OS}. In particular, inequality
(\ref{nbezn0}), (\ref{MC3l2}) hold again, and by lemmas
\ref{short} and \ref{area-kqkqk}, there is a constant $C$
(independent of $M$) such that
\begin{equation}\label{gamma2}
A_2\le C l'(d_3+d_4+1)
\end{equation}
if $l'\ge m/log' n$.

For arbitrary $l'$, by Lemma \ref{area},

\begin{equation}\label{gamma2'}
A_2\le c l'(d_3+d_4+\log'l'),
\end{equation}
and what's more, the logarithmic summand can be replaced by $1$ if
the trapezium $\Gamma_2$ defines a short computation.

To obtain the desired contradiction, we first consider

{\bf Case 1}: Either $l'<n/\log'n$ or the computation
corresponding to $\Gamma_2$ is short.

In view of inequalities (\ref{gamma2'}), the modified (in
comparison with \cite{OS}) task is to show that
\begin{equation}\label{tsel'}
(Mn(n-n_0)+\frac{M}{K^2} l'(l-l'))\ge  C_3l'(d_3+d_4+\log'
l')+C_3(l_3^2+l_4^2)+2\alpha_3 l_3 + 2\alpha_4l_4
\end{equation}
where $C_3$ is a constant that does not depend on $M$, and the
logarithm is replaced by 1 when the trapezium $\Gamma_2$
corresponds to a short computation.

We notice that
\begin{equation}\label{a4}
\frac{M}2 n(n-n_0)\ge C_3 l' log'l'
\end{equation}
if $\Gamma_2$ corresponds to a long computation, because we have
$n-n_0\ge 2$, $n\ge l'\log l'$ and $M\ge C_3$. Then, as in
\cite{OS}, we assume without loss of generality that
$\alpha_3\ge\alpha_4$, and consider two cases.

(a) Suppose we have $\alpha_3\le 2C_3(l-l')$. Then inequalities
(\ref{a122}) and (\ref{a322}) hold as in Lemma \ref{main2}.

The sum of inequalities (\ref{a4}) (if $l'<n/log'n$; otherwise we
do not need it), (\ref{MC3l2}), (\ref{a122}), and (\ref{a322})
gives us both versions of the desired inequality (\ref{tsel'}).

(b) Assume now that $\alpha_3>2C_3(l-l')$. Then to come to a
contradiction, we argue as in case (b) of the proof of Lemma
\ref{main2}, but inequality $\frac{M}2 n(n-n_0)\ge C_3 l' log'l'$
(the analog of (\ref{b322})) follows now just from the assumption
that $l'\le n /log' n$ since $M\ge 2C_3$.

{\bf Case 2}:  $l'\ge n/log'n$ and $\Gamma_2$ corresponds to a
long computation.

Since $n$ has been supposed to be greater than
$10g(g(r-1))\log'n$, we have $l'>10g(g(r-1)$. Then, by Lemma
\ref{area-kqkqk}, we may use inequality (\ref{gamma2}) instead of
(\ref{gamma2'}), which has no term $l'\log l'$. So this term is
absent in (\ref{tsel'}), and we do not need (\ref{a4}), and in
subcase (b), we do not need any analog of inequality (\ref{b322}).

 The lemma is proved by contradiction. \endproof

\section{Proof of Theorem \ref{th1}}

\begin{lemma}\label{combinat}
Let $n$ be the combinatorial perimeter of a reduced diagram
$\Delta$, and $|\partial\Delta|$ the modified perimeter. Then
$n=O(|\partial\Delta|)$ \footnote{We use the Computer Science
``big-O" notation assuming that $f(n)=O(g(n))$ if
$\frac1Cg(n)<f(n)<Cg(n)$ for some positive constant $C$.}
\end{lemma}

\proof As in \cite{OS}, it follows from the definition that
$\delta n\le |\partial\Delta|\le n$. \endproof

{\bf Proof of the theorem.} (1) It is proved in \cite{OS} (Lemma
5.3) that
\begin{equation}\label{quadr}
\eee(\Delta)\le (n/2)^2.
\end{equation}

By lemmas \ref{combinat}, \ref{main2} and inequality
(\ref{quadr}), we have the desired upper bound in the property (1)
of Theorem \ref{th1}. The lower bound follows from the
consideration of the diagrams corresponding to the consequences of
the commutativity relations (\ref{relations}).

(2) We set $n_i=5+4g(i)+2g(g(i))$. Then there is a trapezium
$\Delta$ of height $n_i$ whose area is $O(n_i g(i))=O(n_i\log
n_i)$ and the combinatorial  lengths of top and bottom are equal
to $5+i=O(\log\log n_i)$. ( See Lemma \ref{log} and the comments
to the definition of machine ${\cal M}$ in Section 4.) Lemma
\ref{combinat} and the trick from \cite{OS} with $O(n_i/\log\log
n_i)$ copies of $\Delta$ gluing along sides of these trapezia,
give us a diagram with area $O(\Psi(n_i))$.

(3) Since $g(g(r-1))^2\le g(g(r))$ we conclude from ${\cal
E}(\Delta)\le O(|\partial\Delta|^2)$ and from Lemma \ref{main1}
that $f(n'_i)$ is at most $ O((n'_i)^2)$ for $n'_i=g(g(i))$.

(4) Moreover, the same argument shows that $f(x)$ does not exceed
a quadratic function on the set $\cup_{i=1}^{\infty}
[\frac{d_i}{\lambda_i} ,\; \lambda _i d_i]$, where $d_i =
(n'_i)^{\frac34}$ and $\lambda_i =(n'_i) ^{\varepsilon}$ with
$\varepsilon <1/4$.

(5) It follows from the definitions of $n_i$ and $n'_i$ that
$n_i/3 < n'_i$ for big enough $i$-s. Hence the property (5) of
Theorem \ref{th1} holds with $c_5=1/3$.

Theorem \ref{th1} is proved.

\medskip

{\bf Acknowledgement.} The author is grateful to Mark Sapir for
helpful discussions.

\begin{minipage}[t]{3 in}
\noindent Alexander Yu. Ol'shanskii\\ Department of Mathematics\\
Vanderbilt University \\ alexander.olshanskiy@vanderbilt.edu\\
\end{minipage}
\begin{minipage}[t]{3 in}
\noindent and\\ Department of
Higher Algebra, MEHMAT\\
 Moscow State University\\
olshan@shabol.math.msu.su\\
\end{minipage}

\end{document}